\documentclass{article}

\usepackage{amsmath,amssymb,amsthm}
\usepackage[all]{xypic}

\theoremstyle{definition}
\newtheorem{counter}{}[section]

\newcommand{\nc}{	\newcommand	}
\newcommand{\rc}{	\renewcommand	}

\nc{\ci}{	\circ		}

\nc{\C}{{\Bbb C}}
\nc{\D}{{\Bbb D}}
\nc{\Ee}{{\cal E}}
\nc{\Cc}{{\cal C}}
\nc{\Dd}{{\cal D}}
\nc{\B}{{\cal B}}
\nc{\FF}{{	\cal F	}}
\nc{\G}{{\cal G}}
\nc{\HH}{{\cal H}}
\nc{\N}{{\cal N}}
\nc{\X}{{\cal X}}
\nc{\Ff}{{\cal F}}
\nc{\Aa}{{\cal A}}
\nc{\Oo}{{\cal O}}
\nc{\Pp}{{\cal P}}
\nc{\Ll}{{\cal L}}
\nc{\fg}{{\frak g}}
\nc{\Mm}{{\cal M}}
\nc{\fu}{{	\frak u			}}
\nc{\fr}{{	\frak r			}}
\nc{\fl}{{	\frak l			}}
\nc{\fb}{{	\frak b			}}
\nc{\fn}{{	\frak n			}}
\nc{\fp}{{	\frak p			}}
\nc{\fm}{{	\frak m			}}
\nc{\fh}{{	\frak h			}}
\nc{\ad}{{	\mbox{\rm {Ad}}{}	}}
\nc{\supp}{	\mbox{\rm {supp}}{}      }

\nc{\sub}{	\subset			}
\nc{\al}{	\alpha			}
\nc{\inv}{	^{-1}			}
\nc{\mm}{	\mapsto			}
\nc{\fF}{	\frak F			}
\rc{\k}{\Bbbk}

\begin{document}

\title{Character Sheaves on Reductive Lie Algebras}

\author{I. Mirkovi\'{c}\thanks{Research supported in part by NSF.}}

\date{}

\maketitle

\bigskip

\begin{flushright}
{\large To Borya Feigin}
\end{flushright}

\begin{abstract}

This paper is an introduction,  in a simplified setting, 
to Lusztig's theory of character sheaves.
It develops a notion of character sheaves on 
reductive Lie algebras which is
more general then such notion of Lusztig, and closer to 
Lusztig's theory of character sheaves on groups.
The development is self contained and 
independent of the characteristic $p$ of the ground field.
The results for Lie algebras are then used to
give simple and uniform proofs for some of
Lusztig's results on groups.

\end{abstract}

\section*{Introduction}\

This paper gives a self contained presentation of character sheaves on
reductive Lie algebras, independent of the characteristic $p$ of the field
(except for the section on the characteristic varieties which
have not been defined for $p>0$). For simplicity,
the formulations in  the text involve only the 
$D$-modules (usually irregular),
i.e., the case $p=0$. However
the definitions have ``obvious'' 
modifications 
(in the light of~\cite{Lu3})
for perverse sheaves in characteristic $p>0$ and then the
proofs are the same.

\medskip

The first section studies two Radon transforms and the second defines sheaves
monodromic under an action of a vector group. In the third section we
linearize Lusztig's construction of character sheaves on reductive groups and
show that the resulting character sheaves on Lie algebras are precisely the
Fourier transforms of orbital sheaves. In sections 4 and 5 we give an
elementary proof of results of~\cite{Lu3} (proofs of~\cite{Lu3} are for
sufficiently large $p>0$ and use results from~\cite{Lu1,Lu2}). 
The main result
is that the Fourier transforms of orbital sheaves are precisely the
irreducible constituents
of sheaves obtained by inducing from cuspidal sheaves (5.3). For
$p=0$ character sheaves can be described by having a nilpotent characteristic
variety and a monodromic behavior in some directions (6.4), this is analogous
to a result for groups (\cite{Gi,MV}). 
The last section uses nearby cycles to
reprove some results on induction and restriction from~\cite{Lu2}.

\medskip

The purpose of this paper is to give a  perspective to~\cite{Lu3} and an
introduction to  character sheaves on reductive groups. The original goal
was a direct proof in characteristic zero of the fact
that the cuspidal sheaves on groups
are character sheaves (theorem 6.8). A similar 
proof was found independently by
Ginzburg (\cite{Gi}). Lusztig's original proof is spread though the series of
papers (\cite{Lu2}).

I owe a debt of gratitude to Misha Finkelberg
for his decisive help in bringing this paper to completion.

\subsection*{Notation.}\

We fix an algebraically closed field $\Bbbk$ of characteristic 0. For a
$\Bbbk$-variety $X$,  $\fm (X)$ denotes
the category of holonomic $D$-modules
on $X$ and $D(X)=D^b[\fm(X)]$ is its bounded derived category, our 
basic reference for $D$-modules is~\cite{Bo}. 
If a connected
algebraic group $A$ acts on $X$ then
$\fm_A(X)$ denotes the subcategory of equivariant
sheaves in $\fm (X)$ and $D_A(X)$ is the  corresponding triangulated category
constructed by Bernstein and Lunts (see~\cite{BL,MV}). 
For any morphism of
varieties $X \stackrel{\pi}\longrightarrow Y$ 
there are  direct image
and inverse image functors
$f_*,f_!$
and $f^*,f^!$, and the duality functor
$\D_X:D(X)^o
\to 
D(X)$\
(\cite{Bo}).
If the map $f$ is equivariant under a group $A$ the same functors 
exist on the level of
equivariant derived categories (\cite{BL,MV}).  
We will use a normalized pull-back functor $\pi^{\circ}
\stackrel{\text{def}}= \pi^!$ [$\dim Y-\dim X$], 
and if $\pi$ is an embedding
then $\Ff|X$ will denote the pull-back $\pi^{\circ}\Ff$.
For dual vector bundles
$V$ and $V^*$ one has the Fourier transform equivalence
$\Ff_{V}:
\fm(V)\to \fm(V^*)$, its basic properties
can be found in ~\cite{Br,KL}). 

\medskip

$G$ will denote a reductive algebraic group over $\Bbbk$ and
$P=L\ltimes U$ a Levi decomposition of a parabolic subgroup, while
$\fg,\fp,\fl,\fu$ will be the corresponding Lie algebras. We fix an invariant
non-degenerate bilinear form on $\fg$ and use it to identify $\fg^*$ and
$\fg$, and $\fl^*$ and $\fl$ etc.\footnote{This is just for simplicity of exposition.} 
Now, Fourier transform $\Ff_{\fg}$
is an autoequivalence  of $\fm(\fg)$. 
The adjoint action of  $g \in G$
is denoted $^{g\!} x \stackrel{\text{def}}=(\ad\,g) x$.

\section{
Grothendieck  transform $\G$
and 
horocycle transform $\HH$
}{ }

\bigskip

\counter Let $\B$ be the flag variety of $\fg$, viewed as the moduli of Borel subalgebras of $\fg$. Then $\fb^{\circ}=\{(\fb,x)
\in \B \times \fg,\;x\in \fb\}$ is the 
$G$-homogeneous vector bundle over $\B$ with
the fiber $\fb$ at $\fb \in \B$. 
Similarly, there are vector bundles
$\fn^{\circ}$, $\fg^{\circ}$, $(\fg/\fn)^{\circ}$ with the fibers $[\fb,\fb]$,
$\fg$, $\fg /[\fb,\fb]$ at $\fb \in \B$.

\bigskip

\counter 
Grothendieck resolution of $\fg$ is the projection
$\fb^{\circ}\stackrel{g}\longrightarrow \fg$
while Springer resolution of
the nilpotent cone $\N\subset\fg$ 
is its restriction $\fn^{\circ}=g^{-1}(\N)\stackrel{s}\longrightarrow
\fg$.

\bigskip

\counter Define Grothendieck and horocycle transforms $D(\fb^{\circ})
\stackrel{\G}\longleftarrow D(\fg) \stackrel{\HH}\longrightarrow
D([\fg/\fn]^{\circ})$, 
by $\G=g^!=g^{\circ}$ and $\HH=q_*\circ p^{\circ}$, 
using 
the diagram
\[
\xymatrix@R=1ex{ & \B \times \fg \ar[ddl]_{p} \ar[ddr]^{q} &  & & &(\fb,x)
\ar[ddl] \ar[ddr]& \\
 & & & \mbox{defined by} &  & & \\
 \fg & & (\fg/\fn)^{{\circ}} &  & x & & (\fb,x+[\fb,\fb]).}
\]

Their left adjoints are $\check{\G}=g_!=g_*$ and $\check{\HH}=p_*\circ
q^{\circ}$.

\bigskip

\counter {\bf Theorem.} 
(i) Fourier transform interchanges $\G$ and $\HH$:\
$\FF_{\fb^\ci}\ci \G=\ \HH\ci \FF_\fg$,
hence also $\check{\G}$ and $\check{\HH}$. 

(ii) $\Ff_{\fg}(g_*\Oo_{\fb^{\circ}})=s_*\Oo_{\fn^{\circ}}$.

(iii) $\check{\G}\circ \G = g_* \Oo_{\fb^{\circ}} 
\otimes_{\Oo_{\fg}} -$.

(iv) $\check{\HH}\circ \HH = s_* \Oo_{\fn^{\circ}} * -$.

\bigskip

{\bf Remarks.} a) (ii) is a result of Kashiwara (\cite{Br}) while (iii) (on
$G$ instead of $\fg$) appears in~\cite{Gi,MV}. 
b) The convolution in d) is
defined by $A*B=+_*(A\boxtimes B)$ for the addition $+:\fg\times\fg \to \fg$.

{\it Proof.} 
(i) The map adjoint to 
$\fg^{\circ} 
\stackrel{q}\to
(\fg/\fn)^{\circ}$ is the inclusion 
$\fb^{\circ} 
\stackrel{j}\hookrightarrow
\fg^{\circ}$, so
\begin{align*}
\Ff_\fg\circ\HH 
&=\Ff_{\fg}\circ q_* \circ p^\ci 
=
j^{\circ}\circ\Ff_{\fg^{\circ}}\circ p^{\circ}\\&=j^{\circ}\circ
p^{\circ}\circ\Ff_{\fg}=g^{\circ} \circ \Ff_{\fg}=\G \circ \Ff_{\fg}.
\end{align*}

(ii) 
The sheaves $g_* \Oo_{\fb^{\circ}}$ and $s_* \Oo_{\fn^{\circ}}$ on $\fg$,
are 
direct images from $\fg^{\circ}$ of $(\fb^{\circ}
\stackrel{j}\hookrightarrow \fg^{\circ})_*\Oo_{\fb^{\circ}}$ and
$(\fn^{\circ} \stackrel{j}\hookrightarrow \fg^{\circ})_*\Oo_{\fn^{\circ}}$.
These two are switched by the Fourier transform since $[\fb,\fb]=\fb^{\perp}$
for $\fb \in \B$.

(iii) $A\otimes_{\Oo_{\fg}}g_*
\Oo_{\fb^{\circ}}
=g_*\left(g^{\circ}A\otimes_{\Oo_{\fb^{\circ}}}
\Oo_{\fb^{\circ}}\right)
=\check{\G}(\G A)$.

Now (iv) follows since $\Ff_{\fg}(A*B)
=\Ff_{\fg}A \otimes_{\Oo_{\fg}}
\Ff_{\fg}B$.

\bigskip

\counter {\bf Corollary.} Functors $\check{\G}\circ\G$ and
$\check{\HH}\circ\HH$ both contain identity functors as direct summands.

{\it Proof.} It is well known that $g$ is a small resolution hence $g_*
\Oo_{\fb^{\circ}}$ is a semisimple sheaf and one summand is $\Oo_{\fg}$. Now
the claim for $\G$ follows from (ii) in the theorem and for $\HH$ one uses
(i).

\bigskip

\counter 
Analogous transforms exist for $G$-equivariant sheaves, 
\[
\xymatrix@1{D_G([\fg/\fn]^{\circ}) \ar@<-1ex>[r]_(.585){{\check{\HH}}} & D_G(\fg)
\ar@<1ex>[r]^{{\G}} \ar@<-1ex>[l]_(.4){{\HH}}& D_G(\fb^{\circ})
\ar@<1ex>[l]^{{\check{\G}}}}
;\] 
and from now on we restrict ourselves to the $G$-equivariant setting.

\bigskip

\counter 
A basic tool in the equivariant setting
are
Bernstein's induction
functors $\gamma_B^A,\Gamma_B^A$ (see~\cite{MV}).
If a group $A$ acts on a smooth variety $X$ then for any subgroup
$B$ the forgetful functor 
$\fF^A_B:D_A(X) \to D_B(X)$ has a left adjoint $\gamma_B^A[\dim
A/B]$ which one can describe  via the diagram

\[
\xymatrix{ A \times X \ar[r]^{\nu} \ar[d]_{p} & A \times_B X \ar[d]^{a}\\ X &
X}, \quad \mbox{here} \quad \xymatrix{ (g,x) \ar[r] \ar@{|->}[d] &
\overline{(g,x)} \ar@{|->}[d]\\ x & gx}.
\]
For any $\Aa \in D_B(X)$ there is a unique 
(up to a unique isomorphism)
$\overline{\Aa} \in D_A(A \times_B
X)$ such that $\nu^{\circ}\overline{\Aa}=p^{\circ}\Aa$ in $D_{A\times
B}(A\times X)$. Now $\gamma_B^A \Aa = a_{!}\overline{\Aa}$.
Similarly,
$\Gamma_B^A \Aa = a_{*}\overline{\Aa}$ is a shift of the right adjoint
of the forgetful functor. If $A/B$
is complete then $\gamma_B^A=\Gamma_B^A$.

Since the forgetful functor commutes with all
standard functors
(including the Fourier
transform if $X$ is a vector bundle), the same is also true for $\Gamma_B^A$.

\bigskip

\counter 
Fix a Borel subalgebra $\fb \in \B$, put $\fn=[\fb,\fb]$ and  
consider the maps
\[
\fb^{\circ} \stackrel{r}\hookleftarrow \fb \stackrel{j}\hookrightarrow \fg
\stackrel{\pi}\twoheadrightarrow \fg/\fn \stackrel{s}\hookrightarrow
(\fg/\fn)^{\circ},\quad r(x)=(\fb,x),\;s(y)=(\fb,y). \leqno(*)
\]

\bigskip

\counter {\bf Lemma.} The following diagram is commutative::
\[
\xymatrix@1@R=8ex{ D_G(\fb^{\circ}) \ar@<-1ex>[d]_{r^{{\circ}}}
\ar@<-1ex>[r]_(.5){{\check{\G}}} & D_G(\fg) \ar@<1ex>[r]^(.5){{\HH}}
\ar@<-1ex>[l]_(.5){{\G}} \ar@{=}[d]& D_G((\fg/\fn)^{{\circ}})
\ar@<-1ex>[d]_{s^{{\circ}}} \ar@<1ex>[l]^(.5){{\check{\HH}}} \\ D_B(\fb)
\ar@<-1ex>[u]_{{\Gamma_B^G{\circ}r_*}} \ar@<-1ex>[r]_{{\Gamma_B^G{\circ}j_*}}
& D_G(\fg) \ar@<1ex>[r]^{{\pi_*}} \ar@<-1ex>[l]_{{j^{\circ}}}& D_B(\fg/\fn)
\ar@<-1ex>[u]_{{\Gamma_B^G{\circ}s_*}}
\ar@<1ex>[l]^{{\Gamma_B^G{\circ}\pi^{\circ}}} }.
\]
All vertical arrows are inverse equivalence of categories and
they also commute with the Fourier transform.

{\it Proof.} Since the diagram
$(*)$ is self adjoint, the statements for the first and the second square 
are exchanged by the Fourier transform.
Since $\fb^{\circ}=G\times_B\fb$, the functors $r^{\circ}$ and
$\Gamma_B^G{\circ}r_*$ are inverse equivalences by 
lemma 1.4 in~\cite{MV}. The pull-back
$r^{\circ}$ commutes with the Fourier transform, hence so does
$\Gamma_B^G{\circ}r_*$. For commutativity observe that
$r^{\circ}\circ\G=r^{\circ}\circ g^{\circ}=j^{\circ}$, and 
since $g$ is a $G$-equivariant map
$$
\check{\G}\circ\Gamma_B^G\circ r_*=g_*\circ\Gamma_B^G\circ
r_*=\Gamma_B^G\circ g_*\circ r_*=\Gamma_B^G\circ j_*
.$$

\bigskip

\counter In the remainder we use identifications from Lemma 1.9 as
definitions. So $\fb \stackrel{j}\hookrightarrow \fg \stackrel{\pi}\to
\fg/\fn$ and
\[
\xymatrix@R=1ex{
D_B(\fb) \ar@<-1ex>[ddr]_{\check{\G}
} 
& 
& 
D_B(\fg/\fn)
\ar@<1ex>[ddl]^{\check{\HH}} 
& 
\check{\HH} =\Gamma_B^G\circ\pi^{\circ}, 
& 
\HH= \pi_*\ci \fF^G_B, 
&
	\\ 
&
& 
& 
	\\
& 
D_G(\fg) \ar@<1ex>[uur]^{\HH} \ar@<-1ex>[uul]_{\G} 
& 
& 
\check{\G} 
=\Gamma_B^G\circ j_*, 
&
\G = j^0 \ci \fF^G_B.
}
\]
The theorem 1.4 and its corollary still hold in this setting.

\section{Monodromic sheaves on vector spaces}\

Let $V\subseteq U$ be  finite dimensional vector spaces over $\k$.

\bigskip

\counter Sheaf $0\ne \Aa \in \fm(V)$ is said to be a  character sheaf if 
$+^{\circ}\Aa\cong\Aa \boxtimes \Aa $
for
$+:V\times V \to V$. 
The equivalent condition on
$\B=\Ff_V\Aa $ is $\Delta_*\B 
\cong
\B \boxtimes \B$ for the diagonal $\Delta:V^*
\hookrightarrow V^*\times V^*$. This implies $(\supp\,\B)^2\subseteq
\Delta(V^*) $ hence $\B$ is supported at a point. Now the condition is
equivalent to irreducibility of $\B$. So, the 
character sheaves on $V$ are precisely the
connections $\Ll_{\alpha},\,\alpha \in V^*$,
for
$\Ll_{\alpha}=\Ff_{V^*}[(\alpha \hookrightarrow
V^*)_*\Oo_{\text{pt}}]
$.

\bigskip

\counter Let $U^* \stackrel{q}\to V^*$ be the 
quotient map. For any finite subset
$\theta \subseteq V^*$ we say that a sheaf $\B \in \fm(U)$ is
$\theta$-monodromic if the support of
$\Ff_U\B$ lies in $(-1)\cdot q^{-1}\theta$.
Such sheaves form a Serre subcategory 
$\Mm_{\theta}(U)
$ of 
$\fm(U)$.
The full
triangulated subcategory of $D(U)$ consisting of sheaves with
$\theta$-monodromic cohomologies is precisely the derived category of
$\Mm_{\theta}(U)$.
 Observe that $\Mm_{\theta}(U)
=\bigoplus_{\alpha \in
\theta}\Mm_{\alpha}(U)$. 
For $\alpha\in V^*$, a sheaf $\B \in \fm(U)$ is $\alpha$-monodromic if and
only if $+^{\circ}
\ 
\B\cong\Ll_{\alpha}\boxtimes\B$.

\bigskip

\counter The category of $V$-monodromic sheaves on $U$ is the category
$\Mm(U,V)
\stackrel{\text{def}}
=\bigoplus_{\alpha\in V^*}\Mm_{\alpha}(U)$. 
A sheaf $\B \in
\fm(U)$ is monodromic if and only if the action of 
$V$ on $\Gamma(U,\B)$ 
(via $V\subseteq U\subseteq \Gamma(U,\Dd_U)$), 
is locally finite
--- this is the requirement that the action of 
$V\subseteq\Gamma(U^*,\Oo_{U^*})$
on $\Gamma(U^*,\Ff_U\B)$ is 
locally finite i.e. that the support of $\Ff_U\B$ is
finite modulo $V^{\perp}$. 
So, the category of $V$-monodromic sheaves on $V$ is
semisimple (while this would not be true for a torus).

\section{Character sheaves and the orbital sheaves}\

\counter 
Fix a Borel subalgebra $\fb$ of $\fg$ and denote $[\fn,\fn]$ and
$\fh = \fb/\fn $. Then the diagram
\[
\xymatrix{\fh & \fb \ar
@{->>}[l]_(.4){\nu} \ar@{^{(}->}[r]^(.45){j}& \fg
\ar@{->>}[r]^(.4){\pi} & \fg/\fn \ar@{<-^{)}}[r] & \fh}
\]
is self  adjoint. For any finite 
subset $\theta \subseteq \fh=\fh^*$ let
$\Pp_{\theta}(\fb)
$ be the category of sheaves on $\fb$ supported on
$\nu^{-1}\theta$, so that 
$\Mm_{\theta}(\fg/\fn)
=\Ff_{\fb}(\Pp_{\theta}(\fb))$
is the category of sheaves on $\fg/\fn$
which are $\theta$-monodromic in the direction of $\fh\sub\fg/\fn$.
 Also, we denote
$\Pp(\fb)
=\bigcup_\theta\ 
\Pp_{\theta}(\fb)$, 
$\Mm(\fg/\fn)
=\bigcup_\theta\ 
\Mm_{\theta}(\fg/\fn)
$
and in the equivariant version
$\Mm_{\theta,B}(\fg/\fn)=\Mm_{\theta}(\fg/\fn)\bigcap\fm_B(\fg/\fn)$ etc.

\bigskip

\counter We define character sheaves on $\fg$ as 
irreducible constituents
of the cohomology  sheaves of complexes 
$
\check{\HH}(\Aa)
$ 
for
$\Aa\in\Mm_B(\fb)
=\bigcup_{\theta}\Mm_{\theta,B}(\fb)
$ (see (1.10)). 
On the other hand, an irreducible
$G$-equivariant sheaf on $\fg$ is said to be orbital (\cite{Lu3}), 
if its
support is the closure of a single $G$-orbit.

\bigskip

\counter 
{\bf Lemma} 
(\cite{Lu1}). 
a) Let $x=s+n\in\fg$ be a Jordan
decomposition with $s$ semisimple and $n$ nilpotent. Let $\fp=\fl\ltimes\fu$
be a parabolic subalgebra of $\fg$ such that $\fl=Z_{\fg}(s)$. Then
$^{U\!}x=x+\fu$ for the unipotent subgroup $U$ corresponding to $\fu$.

b) For any $\alpha\in\fh$ 
the
semisimple components of all elements of
$\nu^{-1}\alpha$ are conjugate.

c) For any $G$-orbit $\al$ in $\fg$, 
$\nu(\al\cap \fb)$
is finite.

{\it Proof.} a) appears in the proof of lemma 2.7 in~\cite{Lu1}. For b)
choose semisimple $s\in v^{-1}\alpha=s+\fn$ and put $\fl=Z_{\fg}(s)$ and
$\fp=\fl+\fb=\fl\ltimes\fu$. Now $s$ is the semisimple component of any
element of $s+\fl\cap\fn$ and by a) $(\ad\,U)
(s+\fl\cap\fn)=
(s+\fl\cap\fn)+\fu
=
s+\fn$.

In c) choose $x\in\al$. For $g\in G$,\ 
$^gx\in\fb$ is equivalent to
$x\in 
^{g\inv}\fb$, and then
$\nu(^gx)$ is the same as the image of 
$x$ in
$^{g\inv}\fb/
^{g\inv}\fn\cong \fh$.
So, $\nu(\al\cap\fb)$ is the image of
the fiber
$g\inv x$ under the map $\fb^o\to \fh$ given by
$(\fb',x)\mm x+[\fb',\fb']\in\fb'/[\fb',\fb']\cong \fh$.
This is a complete subvariety of an affine variety, so it is finite.

\bigskip

\counter {\bf Theorem.} 
(See theorem 5.b) in~\cite{Lu3}.) 
Character sheaves
are the same as Fourier transforms of orbital sheaves.

\medskip

\counter {\bf Proposition.} Orbital sheaves are exactly the irreducible
constituents of complexes $\check{\G}(\B)$ with  $\B\in\Pp_B(\fb)$.

{\it Proofs.} The first claim is the Fourier transform of the second since
$\Mm_B(\fg/\fn)
=\Ff_{\fb}\ \Pp_B(\fb)$ 
and 
$\Ff_{\fg}\circ\check{\HH}=
\check{\G}\circ\Ff_{\fg}$ (theorem 1.4). If $\Cc$ is an irreducible
constituent of 
$\check{\G}
(\Dd)$ for some
 $\Dd\in\Pp_B(\fb)$, then this $\Dd$ can be
chosen to be irreducible. Now $\supp(\Dd)\subseteq \nu^{-1} \alpha$ for some
$\alpha\in\fh$. Therefore $\supp(\Cc)$
lies in the closure of $\ad(G)\cdot 
\nu^{-1} \alpha$, and this is
covered by finitely many $G$-orbits (Lemma 3.3.b) and $\Cc$ is orbital.

Conversely, let $\Cc$ be an orbital sheaf. By the Corollary 1.5, $\Cc$ is a
constituent of $\check{\G}(\G(\Cc))$, 
but $\G(\Cc)=j^0\Cc$ is in $\Pp_B(\fb)$ by Lemma 3.3.c.

\bigskip

\counter {\bf Remark.} 
Let us  say that the $L$-packet of character sheaves
attached to an orbit $C$ in $\fg^*$ 
consists of Fourier transforms of 
all orbital
sheaves with the support $\overline{C}$. 
Its elements are in bijection with
the irreducible representations of $\pi_0[Z_G(x)]$ for any $x\in C$.
If $x$ has Jordan decomposition $s+n$ and $L=Z_G(s)$, this is
$\pi_0[Z_L(n)]$.

\medskip

The semisimple part of $C$ is a semisimple orbit in $\fg^*$ so it corresponds
to a Weyl group  orbit $\theta$ in $\fh^*$. We say that
$\theta$ is the
infinitesimal character of sheaves  in the $L$-packet of $C$.
So a character sheaf 
$\Aa$
has infinitesimal character $\theta$ if
$\HH(\Aa)\in\Mm_{\theta,B}(\fg/\fn)$, 
or equivalently, if $\Aa$ is a constituent of some
$\check{\HH}(\B)$ with $\B\in\Mm_{\theta,B}(\fg/\fn)$.

\section{Cuspidal sheaves}\

By observing that the restriction functor commutes with the Fourier transform
we obtain an elementary proof of Lusztig's characterization of cuspidal
sheaves $\Cc$ in terms of the support of $\Cc$ and $\Ff(\Cc)$.

\bigskip

\counter {\bf Restriction and induction.} Let $P$ be a parabolic subgroup of
$G$  with the unipotent radical $U$ and $\overline{P}=P/ U$. Let $\fg$,
$\fp$, $\fu$ and $\overline{\fp}$ be the corresponding Lie algebras.
Restriction (or $\fu$-homology) functor $\mbox{\rm
Res}_{\fp}^{\fg}:D_G(\fg)\to 
D_{\overline{P}}(\overline{\fp})$ is given in terms of the Cartesian diagram 
\[
\xymatrix{ &\fp \ar[dl]_{i} \ar[dr]^{\pi}& & \\ \fg \ar[dr]_{\tau} & &
\overline{\fp}, \ar[dl]^{j} & 
\text{by}\ \ \
{\rm Res}_{\fp}^{\fg}
=
\pi_*
\circ 
i^!\circ \fF^G_P
=j^!\circ\tau_*
\circ \fF^G_P.
\\ & \fg/\fu & &}
\]
For convenience we have chosen to
use the dual of the functor from~\cite{Lu3}, we will see in
(4.7) that this does not matter. 
The left adjoint of
$
{\rm Res}_{\fp}^{\fg}
$
is the induction 
$
\mbox{\rm Ind}_{\fp}^{\fg}
=
\Gamma_P^G\circ i_!\circ\pi^*[\dim G/P]=\Gamma_P^G\circ
i_*\circ\pi^{\circ}$.

\bigskip

\counter {\bf Lemma.} Restriction and induction commute with the Fourier
transform.

{\it Proof.} It suffices to look at the restriction. Since $\tau$ and $j$ are
adjoints of $i$ and $\pi$
\[
\Ff_{\overline{\fp}}\circ(\pi_*\circ i^!\circ \fF^G_P)
=j^{\circ}\circ\tau_*[\dim\fu]\circ \fF^G_P\circ\Ff_{\fg}
=(j^!\circ\tau_*\circ \fF^G_P)\circ\Ff_{\fg}.
\]

The following two observations and their  proofs are from Lemma 2.7
in~\cite{Lu1}.

\bigskip

\counter {\bf Lemma.} Let $x=s+n$ be the Jordan decomposition of $x\in\fg$.
Let $\fp=\fl\ltimes\fu $ be a parabolic subalgebra with a Levi factor
$\ell=Z_{\fg}(s)$. Then for any $\Aa \in\fm(\fg)$
\[
(x\hookrightarrow\fl)^!{\rm Res}_{\fp}^{\fg}\,\Aa=
(x\hookrightarrow\fg)^!\Aa\ [2\dim U]
\]

{\it Proof.} We omit the forgetful functor from the notation and
use the base change
\[
(x\hookrightarrow\fl)^!\mbox{\rm
Res}_{\fp}^{\fg}\,\Aa=(x+\fu\to x)_*(x+\fu \stackrel{i}\hookrightarrow
\fg)^!\Aa
.\] 
Now $x+\fu=\,^{U\!}x$ by (3.3.a) and the $U$-equivariance gives
\[
i^!\Aa=\Oo_{x+\fu}\otimes_{\C}(x\hookrightarrow
x+\fu)^!i^!\Aa=\Oo_{x+\fu}\otimes(x\hookrightarrow\fg)^!\Aa[\dim U]
\]
So the claim follows from
$(\fu\to \text{pt})_*\Oo_{\fu}=\Oo_{\text{pt}}[\dim \fu]$.

\bigskip

\counter {\bf Lemma.} Any sheaf $\Aa\in\fm_G(\fg)$ such that ${\rm
Res}_{\fp}^{\fg}\,\Aa=0$ for all proper parabolics $\fp$, 
is supported in
$Z(\fg)+\N$.

{\it Proof.} For any $0\ne\Aa\in\fm_G(\fg)$ there is a smooth $G$-invariant
subvariety $S\stackrel{j}\hookrightarrow\fg$ open and dense in the support of
$\Aa$ and such that $j^!\Aa$ is a connection $\Ee$ on $S$. For $x$ in $S$
define $s$, $n$, $\fl$, $\fp$ as in (4.3), 
then $(x\hookrightarrow\fl)^!{\rm
Res}_{\fp}^{\fg}\,\Aa= (x\hookrightarrow S )^!\Ee[2\dim \fu]\ne 0$. Therefore
$\fp=\fg$, i.e.,  $s\in Z(\fg)$. So $Z(\fg)+\N$ contains $x$, and then
 also $S$ and
$\overline{S}$.

\bigskip

\counter {\bf Lemma.} Let $\fp=\fl\ltimes\fu$ be a parabolic subalgebra and
$\Aa\in\fm_G(\fg)$. If $\supp(\Aa)\subseteq \N$ then $\supp({\rm
Res}_{\fp}^{\fg}\,\Aa)\subseteq \N\cap\fl$.

{\it Proof.} If $x\in\supp({\rm Res}_{\fp}^{\fg}\,\Aa)\subseteq\fl$ then
$x+\fu$ meets $\supp\,\Aa\subseteq \N$, hence $x \in \N$.

\bigskip

\counter 
An irreducible sheaf $\Aa\in\fm_G(\fg)$ is said to  be cuspidal
(\cite{Lu3}) if (i) ${\rm Res}_{\fp}^{\fg}\,\Aa =0$ for any proper parabolic
subalgebra $\fp$ and (ii) $\Aa=\Ll\boxtimes\B$ for a character sheaf $\Ll$ on
$Z(\fg)$ (see 3.1) and some sheaf $\B$ on $[\fg,\fg]$.

\bigskip

\counter {\bf Theorem.} (see~\cite{Lu3}). Let $\fg$ be a semisimple. An
irreducible sheaf $\Aa\in\fm_G(\fg)$ is cuspidal if and only if $\Aa$ and
$\Ff_{\fg} \Aa$ are supported in $\N$.

{\it Proof.} Since the Fourier transform commutes with the  restriction, if
$\Aa$ is cuspidal so is $\Ff_{\fg} \Aa$. Then they are both supported in $\N$
by lemma 4.4. Conversely suppose now that $\Aa$ and $\Ff_{\fg} \Aa$ live on
$\N$.

For any parabolic $\fp=\fl+\fu$,
both  $\B={\rm Res}_{\fp}^{\fg}\,\Aa$
and $\Ff_{\fl}\B= {\rm Res}_{\fp}^{\fg}(\Ff_{\fg}\Aa)$ are supported in
$\fl\cap\N\subseteq [\fl,\fl]$
by the Lemma 4.5. 
On the other hand if we write
$\B=\overline{\B}\boxtimes(0\hookrightarrow Z(\fl) )_*\Oo_{\text{0}}$ for
a sheaf
$\overline{\B}\in\fm([\fl,\fl])$, 
then $\Ff_{\fl}\B=\Ff_{[\fl,\fl]}\overline{\B}
\boxtimes\Oo_{Z(\fl)}$ has a $Z(\fl)$-invariant support.

Now if $\fp$ is proper then $Z(\fl)\ne 0$, hence $\B=0$. So $\Aa$ is cuspidal.

\bigskip

{\bf Remark.} 
A similar proof was found independently by Ginzburg (~\cite{Gi}).

\bigskip

\counter {\bf Corollary.} Let $\fg$ be semisimple. The set of cuspidal
sheaves is finite and invariant under duality. Any cuspidal sheaf is a
character  sheaf and an orbital sheaf. 
It is $G_m$-monodromic  and has
regular singularities.

{\it Proof.} Since $G$ has finitely many orbits in
$\N$, any  $G$-equivariant sheaf $\B$ supported
on $\N$ is orbital and with regular singularities. 
Since the fundamental groups of $G$-orbits
are finite there are finitely many irreducible $\B$'s. 
The Fourier transform
$\Ff_{\fg}(\Aa)$ 
of a cuspidal sheaf
$\Aa$ 
is supported in $\N$, so it is an orbital sheaf. 
So $\Aa$ is a
character sheaf! Since the duality 
``commutes'' with $\Ff_{\fg}$,\  the dual of $\Aa$
is again cuspidal. Finally, by the $G_m$-invariance of nilpotent orbits
$\Aa$ is smooth on $G_m$-orbits, i.e., $G_m$-monodromic.

\section{Admissible sheaves}

The goal of this section is  to identify two approaches to
character sheaves: 
(i)
by induction from
monodromic sheaves on $\fg/\fn$,
and
(ii) by induction
from  
cuspidal sheaves on Levi factors.
In the standard terminology this is 
the claim that the classes of
{\em character sheaves} and {\em admissible
sheaves}  coincide, and this is essentially  the
theorem 5 from~\cite{Lu3}. 
One direction, that the admissible sheaves are character sheaves
will be essentially obvious 
by now.
The proof of the
converse is based on 
understanding the behavior of character sheaves on
the Lusztig strata in $\fg$, i.e., the behavior
under equisingular change of the semisimple part
of an element of $\fg$.
This is stated as:  
{\em character sheaves} are   {\em quasi-admissible}.
The proof is essentially from~\cite{Lu1}, one
can simplify it here by proving instead the Fourier transform
of the  main result (stated in (5.9)),  
but we use  this version in the next section.

\counter 
Admissible sheaves are defined as irreducible constituents of
all ${\rm Ind}_{\fp}^{\fg}\,\Aa$ for parabolic subalgebras $\fp$ and cuspidal
sheaves $\Aa$ on $\overline{\fp}$ (\cite{Lu3}).

\bigskip

\counter {\bf Lemma.} If $\Aa$ and $\B$ are character (resp. orbital) 
sheaves on
$\fg$ and $\overline{\fp}$, then the 
same holds for all irreducible constituents
of ${\rm Res}_{\fp}^{\fg}\,\Aa$ and ${\rm Ind}_{\fp}^{\fg}\,\B$. The 
complex
${\rm Ind}_{\fp}^{\fg}\,\B$ is semisimple.

{\it Proof.} Fourier transform reduces the lemma to the orbital case.
The semi-simplicity of the induced sheaf
follows from the decomposition theorem since orbital
sheaves are of geometric origin. 
The rest of the proof is the  same as for
proposition 3.5.

\bigskip

\counter {\bf Theorem.} Admissible sheaves are the same as character sheaves.

\bigskip

\counter The class of character sheaves contains 
cuspidal sheaves by (4.8), and is closed under induction 
by (5.2). So it  contains all admissible sheaves.
For the converse we will notice that
character sheaves are quasi-admissible in
(5.6), and use this to show  in (5.8) that character sheaves  are
admissible.

\bigskip

\counter For a Levi subalgebra $\fl$ of $\fg$ we call 
$Z_r(\fl)=\{x\in \fl,\,Z_{\fg}(x)=\fl\}$ the regular part of
the center $Z(\fl)$ of $\fl$. It appears in the subvarieties 
$S_{\fl,\Oo}=\,^{G\!}[Z_r(\fl)+\Oo]$ of $\fg$,
indexed by nilpotent orbits $\Oo$ in
$\fl$.  
This gives the Lusztig  stratification
$\fg=\bigcup\ S_{\fl,\Oo}$ (\cite{Lu1}).

We say that a sheaf $\Aa\in\fm_G(\fg)$ is quasi-admissible if for  any $\fl$
and $\Oo$ all irreducible constituents of
$(Z_r(\fl)+\Oo\hookrightarrow\fg)^!\Aa$ are of the form
$\Ll|Z_r(\fl)\boxtimes\Ee$ for a character sheaf $\Ll$ on $Z(\fl)$ and a
connection $\Ee$ on $\Oo$ (see~\cite{Lu2}).

\bigskip

\counter {\bf Lemma.}  Any character sheaf $\Aa$ is quasi-admissible. In
particular it is smooth on the strata $S_{\fl,\Oo}$.

{\it Proof.} For a pair $\fl,\Oo$ choose a parabolic subalgebra
$\fp=\fl\ltimes\fu$. Now the proof of (4.3) gives
\[
[Z_r(\fl)+\Oo\hookrightarrow\fl]^!\,{\rm Res}_{\fp}^{\fg}\,\Aa
=
[Z_r(\fl)+\Oo\hookrightarrow\fg]^!\Aa\ [2\dim\fu].\leqno{(*)}
\]
By (5.2), irreducible constituents of ${\rm Res}_{\fp}^{\fg}\,\Aa$ are
character sheaves so they are of the form $\Ll\boxtimes\Ff$ for a character sheaf
$\Ll$ on $Z(\fl)$ and $\Ff \in \fm([\fl,\fl])$. Therefore the constituents of
$(*)$ have the needed property.

\bigskip

\counter {\bf Corollary.} Let $\fg$ be semisimple. For an irreducible sheaf
$\Aa\in\fm_G(\fg)$ the following is equivalent: (i) $\Aa$ is cuspidal, (ii)
$\Aa$ is a character sheaf supported on $\N$, (iii) $\Aa$ is both 
a character sheaf
and an orbital sheaf.

{\it Proof.} (i) $\Rightarrow$ (ii) is known and (ii) $\Rightarrow$ 
(iii) is
obvious. Finally, if $\Aa$ is a character sheaf  then lemma 
(5.6) implies that the support of
$\Aa$ is the closure of some stratum $S_{\fl,\Oo}$. If $\Aa$ is also
orbital then $\fl=\fg$ and $\overline{S_{\fl,\Oo}}=\overline{\Oo}\subseteq
\N$. Now (iii) $\Rightarrow$ (i) follows from (4.7) since $\Ff_{\fg}\Aa$ is
also an orbital character sheaf.

\bigskip

\counter {\bf Lemma.} Let $\Aa$ be an irreducible quasi-admissible sheaf on
$\fg$.

(i) $\Aa$ is the irreducible extension of a connection $\Ee$ on one of the
strata $S_{\fl,\Oo}$.

(ii) Connection $\Ee_{\circ}=\Aa|_{Z_r(\fl)+\Oo}$ is irreducible and its
irreducible extension $\Aa_\circ$ to a sheaf on $\fl$ is a constituent of ${\rm
Res}_{\fp}^{\fg}\,\Aa$.

(iii) $\Aa$ is a constituent of ${\rm Ind}_{\fp}^{\fg}(\Aa_{\circ})$.

(iv) If $\Aa$ is also a character sheaf then $\Aa_{\circ}$ is cuspidal.

{\it Proof.} (i) is obvious. 
Denote $S_{\circ}=Z_r(\fl)+\Oo$, then by
$(*)$, we have
$(S_{\circ}\hookrightarrow\fl)^!\,{\rm Res}_{\fp}^{\fg}\,\Aa
=\Ee_{\circ}$
(up
to a shift).  So for (ii) it suffices to see that $S_{\circ}$ is open in
$\supp({\rm Res}_{\fp}^{\fg}\,\Aa)$.

Otherwise there would exist a smooth subvariety $T$ of $\fp-(S_{\circ}+\fu)$
such that $\overline{T}$ meets $S_{\circ}+\fu$ and
$(T\hookrightarrow\fl)^!\,\Aa\ne 0$. One can make $T$ small enough to lie in
some stratum $S_{\tilde{\fl},\tilde{\Oo}}$. Now
$\overline{S_{\tilde{\fl},\tilde{\Oo}}}$ meets
$S_{\circ}+\fu=\,^{U\!}S_{\circ}\subseteq S_{\fl,\Oo}$ and
$(S_{\tilde{\fl},\tilde{\Oo}}\hookrightarrow\fg)^!\,\Aa\ne 0$. This gives
$S_{\tilde{\fl},\tilde{\Oo}}=S_{\fl,\Oo} $.

\medskip

Pick a Cartan subalgebra $\fh_{\circ}$ of $\fl$ and let $W$ and $W_{\fl}$ be
the Weyl groups of $\fg$ and $\fl$. We can suppose that for some $w\in W$ 
subvariety $T$
lies in
\[
^{PwP\!}S_{\circ}=\,^{PwU\!}S_{\circ}=\,^{Pw\!}(S_{\circ}+\fu)
,\]
hence in $\fp\cap
\,^{Pw\!}(S_{\circ}+\fu)=\,^{P\!}[^{w\!}(S_{\circ}+\fu)\cap\fp]=T'$. Let
$A:\fp\to\fh/W_{\fl}$ correspond to $\C[\fh]^{W_{\fl}} \approx
\C[\fp]^{P}\subseteq \C[\fp]$. For any root $\phi$ of $\fh$ in 
$^{w\!}\fl$,
the product
$\psi=\prod_{u\in W_{\fl}}\ u\phi\in \C[\fh/W_{\fl}]$ vanishes on $A(T')$ since the
semisimple parts of elements of 
$^{w\!}(S_{\circ}+\fu)$ lie in
$Z_r(^{w\!}\fl)$. So $\psi\circ A$ vanishes on $T$ and on
the non-empty subvariety $\overline{T}\cap(S_{\circ}+\fu)$.

\medskip

Since for $s\in Z_r(\fl)$ and $n\in\Oo+\fu$ 
one has
$(\psi\circ A
)(s+n)=\phi(s)^{|W_{\fl}|} 
$, 
we find that all roots of $^{w\!}\fl$ vanish at
some $s\in Z_r(\fl)$. Therefore $^{w\!}\fl=\fl$. Now the image of $T'$ in
$\overline{\fp}=\fl$ is $^{L\!}[Z_r(\fl)+\,^{w\!}\Oo]$. Its closure can not
meet $S_{\circ}$ --- otherwise
$^{w\!}\Oo\ne\Oo\subseteq\,\overline{^{w\!}\Oo}$ and $\dim\,
^{w\!}\Oo=\dim\,\Oo$.

(iii) Let $\tilde{S}=G\times_P(S_{\circ}+\fu) \stackrel{\mu}\to S_{\fl,\Oo}$
be the conjugation map. 
Then 
$(S_{\fl,\Oo}\hookrightarrow\fg)^!
\,{\rm Ind}_{\fp}^{\fg}\,\Aa_{\circ}
=\mu_*\mu^{\circ}\Ee
=\Ee\otimes\mu_*\Oo_{\tilde{S}}
$.
Since by lemma 3.3.a map $\mu$ can be identified with 
the map
$G\times_L S_{\circ}\to
G\times_{N_G(L,\Oo)} S_{\circ}$, we see that
$\mu_*\Oo_{\tilde{S}}$ has
$\Oo_{S_{\fl,\Oo}}$ as a summand. 
Finally, since 
$
\supp({\rm Ind}_{\fp}^{\fg}\,\Aa_{\circ})
=\,^{G\!}(\overline{S_{\circ}+\fu})=\overline{S_{\fl,\Oo}}$,
we are done.

(iv) According to  (ii) and (5.2), if $\Aa$ is a character sheaf then so is
$\Aa_{\circ}$. Since $\supp(\Aa_{\circ})\subseteq Z(\fl)$, corollary 5.7
implies that $\Aa_{\circ}$ is cuspidal.

\bigskip

\counter {\bf Corollary.} Orbital sheaves are precisely the 
irreducible constituents
of  sheaves 
${\rm Ind}_{\fp}^{\fg}\,(\Cc\boxtimes \delta_s)$
where  $\Cc$ is a cuspidal sheaf on $[\fl,\fl]$
and $\delta_s$ is the irreducible sheaf supported at a point
$s\in Z(\fl)$.

\section{Characteristic varieties}\

In this section we only consider
$\Bbbk=\C$ since in positive characteristic there is yet no satisfactory
notion of characteristic variety. 
The characteristic variety of $\Ff\in D(\fg)$ is
a subvariety of  $T^*(\fg)\approx \fg\times\fg^*\approx
\fg\times\fg$. We say that
it is nilpotent if it lies in $\fg\times\N$.

\bigskip

\counter Let $\fg$ be semisimple and consider an
irreducible $\Aa\in\fm_G(\fg)$ 
supported in $\N$. Then $\Aa$ is 
$G_m$-monodromic and regular, hence
$\text{pr}_{\fg^*}(\mbox{\rm Ch}\,\Aa)=\supp(\Ff_{\fg}\Aa)$ (\cite{Br}). Therefore
\[
\Aa \quad \mbox{is cuspidal}\quad \Leftrightarrow \quad {\rm Ch} \Aa \quad
\mbox{ is nilpotent.}
\]

\bigskip

\counter 
It is easy to see that the
characteristic variety of any character sheaf is nilpotent. For any
$\B \in \Mm_B$ (section 3.1), 
monodromic property implies that
${\rm Ch}(\B)\subseteq
T^*(\fg/\fn)=\fg/\fn\times\fb$ 
actually lies in $\fg/\fn\times\fn$,
and therefore
${\rm
Ch}[(\fg\to\fg/\fn)^{\circ}\B]\subseteq\fg\times\fn$. So it suffices to 
apply the principle that
${\rm Ch}(\Gamma_B^G\Ff) \subseteq \overline{G\cdot{\rm Ch}(\Ff)}$. 
For $\Ff$ with regular singularities
this is
the lemma 1.2 in~\cite{MV},
however the proof remains correct for irregular $\Dd$-modules 
when the homogeneous space
$G/B$ is complete (the
extension of this lemma to
irregular sheaves and arbitrary  homogeneous spaces
is not true).

\counter {\bf Theorem.} An irreducible sheaf $\Aa\in\fm_G(\fg)$ is a
character sheaf if and only if $\Aa$ is quasi-admissible and ${\rm Ch}(\Aa)$
is nilpotent.

{\it Proof.} It remains to show that a quasi-admissible $\Aa$ with 
${\rm Ch}(\Aa)
\subseteq \fg\times\N$ is a character sheaf. 
Since $\Aa$ is quasi-admissible, we can use  use lemma 5.8
(i)-(iii), its proof and notation. 
It remains to prove a version of
(5.8.iv): if ${\rm Ch}(\Aa)$ is nilpotent then
$\Aa_{\circ}$ is
cuspidal.  This is equivalent (by (6.1)) to:\ 
${\rm Ch}(\Aa_{\circ})$ is
nilpotent.

Since $\Aa$ is quasi-admissible ${\rm Ch}(\Aa_{\circ})$ 
is the union of conormal
bundles to $Z(\fl)+\Oo'$ for some nilpotent
$L$-orbits $\Oo'\subseteq\overline{\Oo}$.
Let $\Oo'$ be one of the orbits contributing to ${\rm Ch}(\Aa_{\circ})$. 
Since
for $z\in Z_r(\fl)$, $n\in\Oo'$ and $x=z+n$, the conormal space at $x$
is
\[
T^*_{S(\fl,\Oo')}(\fg)_x= Z(\fl)^{\perp} \cap
T_x(^{G\!}x)^{\perp}=Z(\fl)^{\perp} \cap Z_{\fg}(x)=Z_{[\fl,\fl]}(n)
=T^*_{Z(\fl)+\Oo'}(\fl)_x
,\]
it suffices to see that $T^*_{S(\fl,\Oo')}(\fg)\subseteq {\rm Ch}(\Aa)$.

Let $\fl\stackrel{j}\hookleftarrow
Z_r(\fl)+\overline{\Oo}\stackrel{i}\hookrightarrow\fg$ 
so that 
$\Aa_{\circ}$ is a
constituent of 
$j_*i^{\circ}\Aa$, hence 
$T^*_{Z(\fl)+\Oo'}(\fl)
\subseteq {\rm Ch}(j_*i^{\circ}\Aa)$. 
To compare $\Aa$ and $i^{\circ}\Aa$ we 
use the radical $U_-$ of the opposite parabolic subgroup
and the map
$\pi :U_-\times U\times Z_r(\fl)\times \overline{\Oo}\to \fg
$,
  defined by
$\pi(\overline{u},u,z,n)=\,^{\overline{u}u\!}(z+n) $. 
The image $S$ of $\pi$  is
open in $\supp(\Aa)$ since
\[
\supp(\Aa)=\overline{S_{\fl,\Oo}}=\,\overline{^{G/P\!}(Z_r(\fl)+\Oo+\fu)}=\,^{G/P\!}
(Z(\fl)+\overline{\Oo}+\fu),
\]
and by lemma 3.3.a\ 
$S=\,^{U_-\!} (Z_r(\fl)+\overline{\Oo}+\fu)$.

The map 
$\pi$ is finite and any $z\in Z_r(\fl)$ has a neighborhood $V$ in
$Z_{r}(\fl)$ such that 
the restriction
of $\pi$ 
to $U_-\times U\times V\times \overline{\Oo}$
is an isomorphism onto its image which is open in $S$. To check this let
$d_i=(\overline{u}_i,u_i,z_i,n_i)\in U_-\times U\times
Z_r(\fl)\times \overline{\Oo}$, $i=1,2,3$. If 
$\pi (d_i)=\pi (d_j)$ then
$s_{ij}=u_i^{-1}\overline{u}_i^{-1}\overline{u}_ju_j\in N_G(L)$ since
$^{s_{ij}\!}z_j=z_i$. 
So if $\pi (d_1)=\pi (d_2)=\pi (d_3)$ and ${s_{12}}L={s_{13}}L$,
then $d_2=d_3$. First,  $s_{23}\in L$, then
$\overline{u}_2^{-1}\overline{u}_3\in U_-\cap P=1$, and so  lemma
3.3.a gives $d_2=d_3 $.

Since $\Aa$ is $G$-equivariant $\pi^{\circ}\Aa=\Oo_{U_-\times
U}\boxtimes i^{\circ}\Aa$. 
Now, the above properties of $\pi$ imply that
$T^*_{S(\fl,\Oo')}(\fg)\subseteq {\rm Ch}(\Aa)$.

\bigskip

\counter For $\Ff\in\fm_G(\fg)$ and $x\in\fg$ 
the fiber
${\rm Ch}_x\Ff
={\rm Ch}\,\Ff\cap T^*_x(\fg)
$
lies in 
the conormal space
$T^*_{G_x}(\fg)_x=Z_{\fg}(x)$, 
so ${\rm Ch}\,\Ff$ is nilpotent if and only if ${\rm Ch}(\Ff)\subseteq \Lambda$ for
$\Lambda=\{(x,y)\in \fg\times\N,\,y\in Z_{\fg}(x)\}$.
The following  is the Lie algebra analogue of Laumon's  result
for groups (\cite{La}).

\bigskip

\counter {\bf Lemma.} $\Lambda$ is a Lagrangian subvariety. It is actually
the union of some of $T^*_{S(\fl,\Oo')}(\fg)$.

\bigskip

{\it Proof.} 
The fiber of $\Lambda$ at $y\in\N$
is $Z_{\fg}(y)=T^*_{G_y}(\fg)_y$. 
So, from the point of view of the projection to $\N$,
$\Lambda$ is the union of  conormal bundles to 
all nilpotent
orbits. So $\Lambda$ is Lagrangian.

Now let $\Oo$ be a nilpotent orbit in a Levi factor $\fl$
and $s\in Z_r(\fl)$, $n\in\Oo$. Then
the tangent space
$T_{s+n}(Z_r(\fl)+\Oo)=Z(\fl)+[n,\fl]$ is orthogonal to
$Z_{\N\cap\fl}(n)=Z_{\N}(s+n)=\Lambda\cap T^*_{s+n}(\fg)$. So
$\Lambda\subseteq \bigcup T^*_{S(\fl,\Oo)}(\fg)$ and since $\Lambda$ is
Lagrangian it is a union of some of $T^*_{S(\fl,\Oo)}(\fg)$.

\bigskip

\counter For sheaves on the group $G$ one defines 
${\rm Res}_P^G$ 
and
cuspidality in an analogous way. By identifying $\N\subseteq\fg$ with the
unipotent cone in $G$ we will not cause any confusion since for sheaves
supported on $\N$, 
the
exponential map intertwines ${\rm Res}_\fp^\fg$ 
and ${\rm Res}_P^G$. We know
that $\supp({\rm Res}_{\fp}^{\fg}\,\Aa)\subseteq\fl\cap\N$ (lemma 4.5) and
similarly on the group. 
In order to calculate ${\rm Res}_{P}^{G}\,\Aa$ (resp. ${\rm
Res}_{\fp}^{\fg}\,\Aa$) we integrate over cosets $e^Y\cdot U$ (resp.
$Y+\fu$), for $Y\in\N\cap\fl$. But $e^Y\cdot U=e^{Y+\fu}$.

\bigskip

\counter {\bf Theorem.} Any cuspidal sheaf $\Aa$ on a semisimple group $G$ has
a nilpotent characteristic variety.
In particular it is a character sheaf.
\bigskip

{\bf Remark.} Lusztig proved in any characteristic that  cuspidal sheaves
on a group are character sheaves (theorem 23.1.b in~\cite{Lu2}). 
The above theorem gives a simpler proof
but only in  characteristic zero.
A similar proof was found
by Ginzburg (\cite{Gi}).

{\it Proof.} 
The second sentence follows from the first 
since any
irreducible $G$-equivariant regular $D$-module on $G$ with a nilpotent
characteristic variety is a character sheaf (\cite{Gi,MV}).
The first claim is that
${\rm Ch}(\Aa)\subseteq T^*(G)=G\times\fg^*=G\times\fg$ 
actually lies in $G\times\N$. If 
$\supp(\Aa)\sub\N $ then the claim follows
from (6.7) and (6.1). To reduce the situation to this case 
we follow~\cite{Lu1}.

The pull-back of $\Aa$ to any cover of $G$ is obviously cuspidal so we can
suppose that $G$ is simply connected. Choose $S$, $\Ee$ and $x=sn\in S$ as in
the proof of Lemma 4.4. Then $H=Z_G(s)$ is connected. It cannot lie in a
proper Levi subgroup $L$ of $G$ since we could repeat the proof of (4.4) and
get ${\rm Res}_{P}^{G}\,\Aa\ne 0$. Therefore $H$ is semisimple and $s$ 
is of
finite order. So $G$ has finitely many orbits in
$S$ and we can assume
that $S=\,^{G\!}x$ is a single orbit.

Let $\Oo=\,^{H\!}n$ and let $\Aa_H$ be the irreducible extension of the
connection $(s^{-1})^{\circ}(\Ee|s\Oo)$ from $\Oo$ to $H$. Since\
$\supp(\Aa)=\overline{S}=\,^{G\!}(s\overline{\Oo})\approx G\times_H
\overline{\Oo}$, we see that 
${\rm Ch}(\Aa)$ is the union of 
conormal bundles for subvarieties $^G(s{\Oo}')\sub G$ 
over all 
 $H$-orbits $\Oo'$  in
$\overline{\Oo}$ 
such that 
$
T^*_{s\Oo'}(\fl)\subseteq {\rm Ch}(\Aa_H)
$. 
For any
$v\in\Oo'$
\[
T^*_{^G(s\Oo')}(G)_{sv}=Z_{\fg}(sv)=Z_{\fh}(v)=T^*_{\Oo'}(H)_v,
\]
so it remains to see that $\Aa_H$ is cuspidal.

Let $P=L\ltimes U$ be a parabolic subgroup of $G$ such that $P_H=L_H\ltimes
U_H$ (for $P_H=P\cap H$ etc.) is a parabolic subgroup of $H$. Let
$v\in\overline{\Oo}\cap P_H$, $T_H=vU_H\cap\overline{\Oo}$ and
$T=svU\cap\,^{G\!}(s\overline{\Oo})$. The projection to the semisimple part of the Jordan decomposition gives an algebraic map
$T\stackrel{\alpha}\to\,^{G\!}s\cap P$ 
(it is a restriction of the map
$\,^{G\!}(s\overline{\Oo})=G\times_H s\overline{\Oo}\to G\times_H
s=\,^{G\!}s$). By composing $\alpha$ with $P\to \overline{P}$ we find that
${\rm Im}(\alpha)\subseteq \,^{G\!}s\cap sU$. Since $\,^{U\!}s$ is a
connected component of $\,^{G\!}s\cap sU$
\[
\alpha^{-1}(^{U\!}s)=\,^{U\!}(\alpha^{-1}s)=\,^{U\!}(svU\cap
s\overline{\Oo})=\,^{U\!}(sT_H)\approx U\times_{U_H}sT_H
\]
is open in $T$. Therefore the integral of $\Aa$ over $\alpha^{-1}(^{U\!}s)$
vanishes. By $U$-equivariance the same is true for the integral $\Aa$ over
$sT_H$. 
Finally, since $(sT_H\hookrightarrow G)^!\,\Aa=
(T_H\hookrightarrow H)^!\,\Aa_H$
(up to a shift), $\Aa_H$ is cuspidal.

\section{An application of nearby cycles}{ }

Let $C$ be a smooth curve and $O\in C$. 
Any function $f:X\to C$ defines exact functors
of nearby and vanishing
cycles 
$\fm(X)\stackrel{\psi_f,\phi_f}\longrightarrow
\fm_{f^{-1}(0)}(X)$, from $D$-modules on $X$ to the subcategory of
$D$-modules on $X$ supported in the fiber 
$f^{-1}(O)$ 
(\cite{Be}). Suppose that a group $G$ acts on $X$ and fixes the function $f$.

\bigskip

\counter {\bf Lemma.} For any subgroup
$B$ such that  $G/B$ is complete, functor $\Gamma_B^G$ 
commutes with
the equivariant versions of
$\psi_f$ and $\phi_f$.

{\it Proof.} 
Computing $\Gamma_B^G$ involves inverse images $\alpha^{\circ}$
for smooth maps $\alpha$ and direct images $\beta_*$ for proper maps $\beta$
(see (1.7)), and these commute with $\psi_f$ and $\phi_f$.

\bigskip

\counter As a consequence nearby and vanishing cycles
will commute with $\check{\HH}$,
$\check{\G} $ and induction. This gives a simple proof for the following
result of Lusztig.

\bigskip

\counter {\bf Theorem.} (\cite{Lu2}). Let $\Aa$ and $\B$ be 
character sheaves
(resp. orbital sheaves) on $\fg$ and $\fl$. Then

(i) ${\rm Ind}_{\fp}^{\fg}\,\B$ is a semisimple sheaf,

(ii) ${\rm Res}_{\fp}^{\fg}\,\Aa\in D_{\geqslant 0}(\fl)$,

(iii) ${\rm Ind}_{\fp}^{\fg}\,\B$ and ${\rm Res}_{\fp}^{\fg}\,\Aa$ do not
depend on the choice of a parabolic subalgebra $\fp=\fl\ltimes\fu$.

{\it Proof.} As usual it suffices to look at the orbital sheaves. 
Since $\B$
is an orbital sheaf on $\fl$ we have
$\B=\delta_s\boxtimes\Cc$ for an orbital sheaf
$\Cc\in\fm_L([\fl,\fl])$ and
$\delta_s=(s\hookrightarrow\fl)_*\,\Oo_{\text{pt}}$
for some
$s\in Z(\fl)$. 
Suppose that
$\Cc$ is supported in
$\N\cap\fl$ 
and  $s\in Z_r(\fl)$.
Then $^{G\!}(s+\N\cap\fl)=G\times L (s+\N\cap\fl)$ and
for any
$n\in\N\cap\fl$ we have
$s+n+\fu=\,^{U\!}(s+n)$. 
This
shows that for $\xymatrix{\fl \ar@<-0.7ex>[r]_{j}&\fp
\ar@<-0.7ex>[l]_{p} \ar@{^(->}[r]^{i} &\fg}$
\[
{\rm Ind}_{\fp}^{\fg}\,\B
=
\Gamma_P^Gi_*(p^{\circ}\B)
=
\Gamma_P^Gi_*(\Gamma_L^Pj_*\B)
=
\Gamma_L^G(ij)_*\B
\]
is a sheaf independent of $\fp$. The point was that in this case 
the parabolic induction is just the naive induction from
a Levi factor:\
${\rm Ind}_{\fp}^{\fg}\,\B={\rm Ind}_{\fl}^{\fg}\,\B$.

For the general $s\in Z(\fl)$ we construct a sheaf on 
$\fl\times{\Bbb A}^1$ by
$\tilde{\B}=\Cc\boxtimes\delta_K$ for $K=\{(s+c\alpha,c),\,c\in k\}\subseteq
Z(\fl)\times {\Bbb A}^1$ such that 
$s+c\alpha\in Z_r(\fl)$
for $c$ in some open dense $U\subseteq
{\Bbb A}^1$. By the above calculation,  sheaf ${\rm Ind}\,\tilde{\B}$ is
independent of $\fp$ on $\fg\times
U$. Hence so is
$\psi({\rm Ind}\,\tilde{\B})={\rm Ind}(\psi\tilde{\B})={\rm Ind}\,\B$.

By the transitivity of induction and by corollary (5.9), 
claim (i) and the induction
part of (iii) follow for all orbital sheaves. 
Since ${\rm Ind}_{\fp}^{\fg}$ and
${\rm Res}_{\fp}^{\fg}$ are adjoint functors between full subcategories of
$D_G(\fg)$ and $D_L(\fl)$ consisting of complexes supported on finitely many
orbits we also get the restriction part of (iii). 
The claim  (ii)
also follows by adjunction since for $i<0$ we have
$0={\rm
Hom}\left({\rm Ind}_{\fp}^{\fg}\,\B,\,\Aa[i]\right)={\rm Hom}\left(\B,\,{\rm
Res}_{\fp}^{\fg}\,\Aa[i]\right)$, hence
$H^j({\rm
Res}_{\fp}^{\fg}\,\Aa)=0$ for $i<0$.

\bigskip

\counter {\bf Remarks.} (i) The proof is self contained for the class of
sheaves supported in nilpotent cones. This can be used for a proof of the
Theorem 5.3 which avoids Lemma 5.8. (ii) Actually ${\rm Res}_{\fp}
^{\fg}\,\Aa$ is also a semi-simple sheaf (\cite{Lu2}). 
A simple proof in characteristic zero was found by Ginzburg (\cite{Gi1}).

\bigskip

\counter Let $s:\check{\N}\to\fg$ be the Springer resolution of $\N$. Then
the proof above shows that the Springer sheaf
$s_*\Oo_{\check{\N}}$ is the limit
$\lim\limits_{C\to 0}\delta_C$ of  
delta-distributions
$\delta_C=(C\hookrightarrow\fg)_*\Oo_C$ 
 on regular
semisimple conjugacy classes $C$.
More precisely, 
 any 
regular
semisimple conjugacy class
$C$ defines  a $G_m$ family 
of $D$-modules 
$\delta_{s\cdot C},\ s\in G_m$,
and the nearby cycle limit of the family at $s=0$ is the Springer sheaf.

\noindent
{\small 
DEPT. of MATHEMATICS, UNIVERSITY of MASSACHUSETTS,
AMHERST,  MA 01003 USA}
\newline
{\small		
Email address: mirkovic@math.umass.edu	}

\end{document}